\documentclass{amsart}
\usepackage{amsfonts,amsmath}
\usepackage{graphicx,rotate}
\usepackage{epic,eepic}

\newtheorem{theorem}{Theorem}
\newtheorem{claim}{Claim}

\def\G{{\mathcal G}}
\def\F{{\mathcal F}}

\title[The Tits alternative for non-spherical Pride groups]
{The Tits alternative for non-spherical Pride groups} % This is the full title of the paper

\author{Natalia Kopteva} 
\address{Sobolev Institute of Mathematics, Acad. Koptyug ave., 4,
 Novosibirsk 630090, RUSSIA}
\email{natasha@math.nsc.ru}

\author{Gerald Williams}
\address{Institute of Mathematics, Statistics, and
Actuarial Science, University of Kent, Canterbury, Kent CT2 7NF, UK}
\email{g.williams@kent.ac.uk}

\thanks{The first author was partially supported by FP6 Marie
Curie IIF Fellowship.}

\date{\today}



\begin{document}

\maketitle

\begin{abstract}
Pride groups, or ``groups given by presentations in which each
defining relator involves at most two types of generators'', include
Coxeter groups, Artin groups, triangles of groups, and Vinberg's
groups defined by periodic paired relations. We show that every
non-spherical Pride group that is not a triangle of groups satisfies
the Tits  alternative.
\end{abstract}

\section{Introduction}

Pride groups, or ``groups given by presentations in which each defining
relator involves at most two types of generators'' \cite{Pri92},
include Coxeter groups, Artin groups, triangles of groups, and
Vinberg's groups defined by periodic paired relations.
The cohomology of Pride groups was considered in \cite{Pri92},
geometric invariants were considered in \cite{Mei97}, and
a Freiheitssatz was proven in \cite{Cor96}.

In this paper we consider the Tits alternative for the class of Pride groups.
Recall that a class of groups $\mathcal C$  satisfies the
{\it Tits alternative} if each group in $\mathcal C$ contains a
non-abelian free subgroup or has a soluble subgroup of finite index.
This property is named after Tits who established
that it is satisfied by the class of linear groups~\cite{Tit72};
in particular, it holds for Coxeter groups.

The Tits alternative has been considered, for example, for the
classes of mapping class groups of compact
surfaces~\cite{Iva84,McC85}, outer automorphism groups of free
groups of finite rank~\cite{BFHa,BFHb}, subgroups of Gromov
hyperbolic groups~\cite{Gro87}, groups acting on ${\rm CAT}(0)$
cubical complexes~\cite{SW05}, triangles of groups~\cite{HK06}, and
groups defined by periodic paired relations~\cite{Vin97,Wil06}.

In this paper we prove the following

\begin{theorem}\label{thm}
Every non-spherical Pride group $G$ based on a
graph with at least 4
vertices contains a non-abelian free subgroup, unless it
is based on the graph shown in Figure~\ref{square}, in which
case $G$ is virtually abelian and has presentation
$$
\langle x_1,x_2,x_3,x_4\, |\, x_1^2,x_2^2,x_3^2,x_4^2,(x_1 x_2)^2,
(x_2 x_3)^2,(x_3 x_4)^2,(x_4 x_1)^2\rangle.
$$
\end{theorem}

It is interesting to note that the ``negatively curved'' property of
containing a non-abelian free subgroup is found in this non-positively
curved class of groups.

\begin{figure}[htbp]
\centering
\tiny
\setlength{\unitlength}{0.2mm}
\begin{picture}(200,200)(10,-15)
\drawline(10,0)(210,0)
\drawline(10,0)(110,171)
\drawline(110,171)(210,0)
\drawline(10,0)(110,57)
\drawline(110,171)(110,57)
\drawline(110,57)(210,0)
\put(10,0){\circle*{8}}
\put(210,0){\circle*{8}}
\put(110,57){\circle*{8}}
\put(110,171){\circle*{8}}
\put(115,63){\makebox(0,0)[lb]{\smash{$\langle x_1\,|\,x_1^2\rangle$}}}
\put(-20,-15){\makebox(0,0)[lb]{\smash{$\langle x_3\,|\,x_3^2\rangle$}}}
\put(200,-15){\makebox(0,0)[lb]{\smash{$\langle x_2\,|\,x_2^2\rangle$}}}
\put(95,180){\makebox(0,0)[lb]{\smash{$\langle x_4\,|\,x_4^2\rangle$}}}
\put(95,20){\makebox(0,0)[lb]{\smash{$\{(x_1x_2)^2\}$}}}
\put(80,-15){\makebox(0,0)[lb]{\smash{$\{(x_2x_3)^2\}$}}}
\put(0,90){\makebox(0,0)[lb]{\smash{$\{(x_3x_4)^2\}$}}}
\put(92,75){\makebox(0,0)[lb]{\smash{\begin{rotate}{$\{(x_4x_1)^2\}$}\end{rotate}}}}
\put(60,35){\makebox(0,0)[lb]{\smash{$\emptyset$}}}
\put(165,90){\makebox(0,0)[lb]{\smash{$\emptyset$}}}
\end{picture}
\normalsize
\caption{}\label{square}
\end{figure}

We now give our formal definitions.
Let $\G$ be a finite simplicial graph with vertex set
$I(\G)$, and edge set $E(\G)$.
Further, let there be non-trivial groups
$G_i$ (with fixed finite presentations) associated
to each vertex $i\in I(\G)$ and, in addition, for each edge
$\{i,j\}\in E(\G)$ let $R_{ij}$ be a (possibly empty)
finite collection of cyclically reduced words. We assume each word in
 $R_{ij}$ is of free product length greater than or equal to 2
 in $G_i*G_j$.
The {\it Pride group}
based on the graph $\G$ with groups $G_i$ assigned to the vertices
and with edge relations
$R=\cup_{\{i,j\}\in E(\G)} R_{ij}$
is the group $G=*_{i\in I(\G)} G_i/N$, where $N$ is the normal
closure of $R$
in $*_{i\in I(\G)} G_i$.

We refer to the groups $G_i$ as {\it vertex groups},
and we define the {\it edge groups} to
be $G_{ij}= \{G_i * G_j\}/N_{ij}$,
where $\{i,j\}\in E(\G)$ and where $N_{ij}$ is the normal closure
of $R_{ij}$ in $G_i* G_j$. More generally,
if $\F$ is any full subgraph of $\G$ with vertex set
$I(\F)\subseteq I(\G)$, then the
{\it subgraph group} $G_{\F}$ is
$\{*_{i\in I(\F)} G_i\}/\{N_{ij}|\{i,j\}\in E(\F)\}$.
In particular, $G_{\G} = G$.

For each $i,j\in I(\G)$, the natural homomorphisms
$G_i\to G_{ij}$, $G_j\to G_{ij}$
determine a homomorphism $G_i * G_j \to G_{ij}$.
Let $m_{ij}$ denote the length of a shortest
non-trivial element in its kernel
(in the usual length function on the free product), or
put $m_{ij} =\infty$ if the kernel is trivial.
Note that either $m_{ij} = 1$ (in which case one of
the natural maps $G_i\to G_{ij}$, $G_j\to G_{ij}$
is not injective), or $m_{ij}$ is even or infinite.
The {\it Gersten-Stallings angle} $(G_{ij};G_i,G_j)$
between the groups $G_i$ and $G_j$ in the group $G_{ij}$
is defined to be $2\pi/m_{ij}$ for $m_{ij}>1$, and 0 for
$m_{ij}=\infty$ \cite{Sta91}.

In \cite{Pri92} Pride formulated the following
asphericity condition.
A Pride group $G$ based on a graph $\G$
(with $|I(\mathcal G)|\geq 3$) 
is said to be {\it non-spherical} if
\begin{itemize}
\item[(i)]
$(G_{ij};G_i,G_j)\leq\pi/2$ for all $i,j\in I(\G)$; and
\item[(ii)]
for any triangle $\{i,j,k\}$ in $\G$
$$
(G_{ij};G_i,G_j)+(G_{jk};G_j,G_k)+(G_{ik};G_i,G_k)\leq \pi.
$$
\end{itemize}

In the non-spherical case we can assume that the graph $\G$ is
complete. To see this, observe that if $i,j\in I(\G)$ and
$\{i,j\}\notin E(\G)$ then we can add the edge $\{i,j\}$
and set $R_{ij}=\emptyset$ without changing the group~$G$.

If $|I(\G)|=3$ then the Pride group $G$ is the colimit
of a triangle of groups.
In \cite{HK06}, it was proved that if the angle sum of the 
triangle is strictly less than $\pi$ then $G$
contains a non-abelian free subgroup. 
In the same paper the Tits alternative was proved 
for a particular class of non-spherical triangles of groups,
 namely,
for non-spherical generalized tetrahedron groups. In general, 
it is unknown if this property holds for non-spherical 
triangles of groups.

We also remark that every Pride group in which $m_{ij}>1$ for all $i,j$ can
be represented in terms of a 2-complex of groups. Moreover, if the Pride
group is non-spherical then the corresponding complex can be chosen to
be non-spherical.

\section{Proof of Theorem~\ref{thm}}

Our method of proof has evolved from that developed
in \cite{ERST} and~\cite{HK06}.

Let $G=G_{\G}$ be a non-spherical Pride group, where $\G$ is complete.
First suppose that $\G$ has four vertices.
Let $I(\G)=\{1,2,3,4\}$ and
let $X=G_1$, $Y=G_2$, $Z=G_3$ and $T=G_4$.
We shall sometimes write $G_{XY}$ for $G_{12}$, $G_{XZ}$ for $G_{13}$
and so on.
Label the vertices of $\G$ by the vertex groups and each edge $\{i,j\}$
by $(G_{ij};G_i,G_j)$.

If $(G_{ij};G_i,G_j)+(G_{jk};G_j,G_k)+(G_{ik};G_i,G_k)< \pi$
for some $\{i,j,k\}\subset I(\G)$ then, by~\cite{HK06},
$G_{ijk}$ contains a non-abelian free subgroup.
By \cite{Cor96}, every subgraph group embeds, so
$G$ also contains a non-abelian free subgroup.
Hence, we may assume that for all $i,j,k\in I(\G)$
the angle sum is exactly~$\pi$. 

Suppose that the edges incident to $T$ are labelled by
 $\theta$, $\alpha$, and $\beta$.
Since the angle sum is $\pi$ for each triangle 
it follows that the edges that do not
share any vertices have the same labels
and all triangles in $\G$ are labelled by one of
$\{\theta,\alpha,\beta\}=\{\pi/2,\pi/2,0\}$,
$\{\pi/2,\pi/3,\pi/6\}$, $\{\pi/2,\pi/4,\pi/4\}$,
$\{\pi/3,\pi/3,\pi/3\}$.
Without loss of generality we may assume that
$\theta\geq \alpha\geq \beta$ and that
\begin{eqnarray*}
&&(G_{XZ};G_X,G_Z)=(G_{YT};G_Y,G_T)=\theta,\\
&&(G_{XY};G_X,G_Y)=(G_{ZT};G_Z,G_T)=\alpha,\\
&&(G_{YZ};G_Y,G_Z)=(G_{XT};G_X,G_T)=\beta.
\end{eqnarray*}

Suppose $(\theta,\alpha,\beta)\not=(\pi/2,\pi/2,0)$ and
consider a presentation $\mathcal{P}$
for $G$.
Since all the vertex groups are non-trivial, we may choose non-trivial
elements $x\in X$, $y\in Y$, $z\in Z$ and $t\in T$
such that $x,y,z,t$ are all generators of $\mathcal{P}$.
We shall show that $u=xyztxyz$ has infinite order in $G$ and
that $t$ and $u$ generate a free product.

Let $w(t,u)=t^{p_1}u^{q_1}\dots t^{p_m}u^{q_m}$ or $w(t,u)=u^{q_1}$,
where $m\geq 1$ and each $p_i,q_i\not=0$,
and assume that  $w(t,u)=1$ in~$G$.
Consider a van Kampen diagram $K$ over $\mathcal{P}$
whose boundary label
is $w(t,u)$.
Let $D$ be an extremal disk of $K$.
We divide $D$ into $G_{ij}$-regions.
If two $G_{ij}$-regions intersect at least at one edge, then we can
amalgamate them into a single region. We continue in this
way as often as possible, and so get a division of
$D$ into {\it maximal $G_{ij}$-regions}.
(Note that the resulting division of $D$ is not necessarily unique.)

By \cite{Cor96}, the edge groups embed, so it
can be assumed that the maximal regions are simply connected.
Let $\widehat D$ be the resulting diagram.
On the boundary of $\widehat D$
an edge of $\widehat D$ is defined to be
a longest path whose edges are labelled by elements of the
same vertex group.
In the interior an edge is defined to be
the intersection of two adjacent maximal $G_{ij}$-
and $G_{ik}$-regions. Note that it is a path labelled by
elements of~$G_i$.

Now place $\widehat D$ on the sphere and take its dual $D^*$.
Let $v_0$ be the vertex corresponding to
${\mathbb S}^2\backslash \widehat D$.
We call a region of $D^*$ {\it exterior} if it involves $v_0$
and {\it interior} otherwise.
We give each corner at a vertex of $D^*$ of degree $\delta$ the angle
$2\pi/\delta$. The curvature $c(\Delta)$ of a region $\Delta$ of
degree $q$ whose vertices have degrees 
$\delta_1, \delta_2, \dots, \delta_q$
is then defined by
$$
c(\Delta)=(2-q)\pi+\sum_{i=1}^q \frac{2\pi}{\delta_i}.
$$
Then
$$\sum_{\Delta\subset D^*}c(\Delta)=2\pi\chi({\mathbb S}^2)=4\pi.$$

We remark that one can use the Gersten-Stallings angles to
estimate the curvature as follows.
Suppose that a vertex $v\not=v_0$ of $D^*$ 
has degree $\delta$ and
comes from a maximal $G_{ij}$-region of~$\widehat D$.
Since the Gersten-Stallings angle $(G_{ij};G_i,G_j)$
is defined to be $2\pi/m_{ij}$, where $m_{ij}$ is
the length of a shortest non-trivial element in the kernel
of $G_i*G_j\to G_{ij}$, 
we have that $\delta\geq m_{ij}$ so 
$2\pi/\delta \leq (G_{ij};G_i,G_j)$. Moreover,
the non-spherical condition implies 
$(G_{ij};G_i,G_j)\leq \pi/2$ so $\delta \geq 4$.

Let $\Delta$ be an interior region of $D^*$ of degree $q$.
Observe that $q\geq 3$.
If $q\geq 4$ then
$$
c(\Delta)\leq (2-q)\pi+\sum_{i=1}^q \frac{\pi}{2}\leq 0
$$
and if $q=3$ then (for some distinct $i,j,k\in I(\G)$)
$$
c(\Delta)\leq
-\pi+(G_{ij};G_i,G_j)+(G_{jk};G_j,G_k)+(G_{ik};G_i,G_k)=0.
$$
Thus, the sum of the curvatures of interior regions is non-positive.

Consider exterior regions.
Observe that an exterior region can be a 2-gon.
It is convenient to define
$d(\Delta)=c(\Delta)-2\pi/N$, where $N={\rm deg}(v_0)$
is the number of exterior regions. Then
$$
\sum_{\mathrm{ext. }\Delta}c(\Delta)=
\sum_{\mathrm{ext. }\Delta}d(\Delta) +2\pi.
$$
We aim to show $\sum_{\mathrm{ext. }\Delta}d(\Delta)\leq 0$
and obtain a contradiction.

Split the boundary of $\widehat D$ into
$t^{p_i}$ and
$u^\varepsilon=(xyztxyz)^\varepsilon$ pieces, where $\varepsilon=\pm1$.
We now consider the sum of the curvatures
of the exterior regions of $D^*$ arising from each
$u^\varepsilon$ piece.

\begin{figure}[htbp]
\centering
\begin{tabular}{cc}
\setlength{\unitlength}{0.15mm}
\begin{picture}(200,200)(10,-15)
\drawline(10,0)(210,0)
\drawline(10,0)(110,171)
\drawline(110,171)(210,0)
\drawline(10,0)(110,57)
\drawline(110,171)(110,57)
\drawline(110,57)(210,0)
\put(10,0){\circle*{8}}
\put(210,0){\circle*{8}}
\put(110,57){\circle*{8}}
\put(110,171){\circle*{8}}
\put(90,63){\makebox(0,0)[lb]{\smash{$T$}}}
\put(-15,-18){\makebox(0,0)[lb]{\smash{$A$}}}
\put(215,-18){\makebox(0,0)[lb]{\smash{$Y$}}}
\put(115,175){\makebox(0,0)[lb]{\smash{$B$}}}
\put(65,20){\makebox(0,0)[lb]{\smash{$\psi$}}}
\put(100,-18){\makebox(0,0)[lb]{\smash{$\phi$}}}
\put(165,90){\makebox(0,0)[lb]{\smash{$\psi$}}}
\put(115,100){\makebox(0,0)[lb]{\smash{$\phi$}}}
\put(130,18){\makebox(0,0)[lb]{\smash{$\theta$}}}
\put(45,90){\makebox(0,0)[lb]{\smash{$\theta$}}}
\end{picture}
\quad & \quad
\includegraphics[width=7 cm]{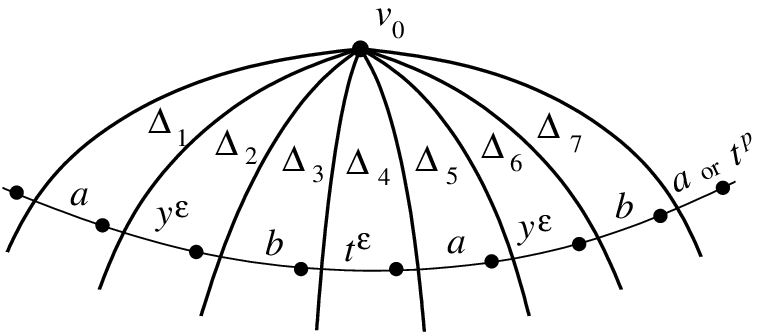}
\\
(a) \quad & \quad (b)
\end{tabular}
\caption{}\label{upiece}
\end{figure}

The analysis of a $u^\varepsilon$ piece is essentially the same
for both $\varepsilon=+1$ and $\varepsilon=-1$, so we
introduce the following notation.
Let $(a,b)=(x,z)$, $(A,B)=(X,Z)$, $(\phi,\psi)=(\alpha,\beta)$
if $\varepsilon=+1$
and let $(a,b)=(z^{-1},x^{-1})$, $(A,B)=(Z,X)$,
$(\phi,\psi)=(\beta,\alpha)$ if $\varepsilon=-1$.
Note that with this convention $\theta\geq\phi$ and $\theta\geq\psi$.
Figure~\ref{upiece}(a) indicates the Gersten-Stallings angles
between the vertex groups and
Figure~\ref{upiece}(b) shows the form of a $u^\varepsilon$ piece,
where the exterior regions of
$D^*$ are labelled by $\Delta_i$ ($1\leq i\leq 7$) and each
$\Delta_i$ is a $q_i$-gon.

Since any three consecutive edges on the boundary of
$\widehat D$ are labelled by elements of three different
vertex groups, no two exterior 2-gons of $D^*$ can be adjacent. 
Therefore, at most four of the $\Delta_i$ can be 2-gons.
Denote the chain $\Delta_1\Delta_2\dots\Delta_7$ by $S$ and
write $d(S)=\sum_{i=1}^7 d(\Delta_i)$.
Denote by $v_1$ the vertex of $\Delta_1\backslash\Delta_2$
adjacent to~$v_0$.

We shall make frequent use of the following observations.
Let $\Delta$ be an exterior $q$-gon. If $q=2$ then
$d(\Delta)=2\pi/\delta\leq\pi/2$. 
If $q=3$ then no two adjacent vertices
of~$\Delta$ arise from maximal $G_{AB}$- or $G_{YT}$-regions, and so
$d(\Delta)\leq-\pi+\pi/2+\pi/3=-\pi/6$. Similarly, if $q\geq 4$ then
$d(\Delta)\leq -2\pi/3$.

\begin{claim}\label{claim1}
If $v_1$ does not arise from any maximal $G_{AT}$-region
then $d(S)\leq 0$.
\end{claim}

\begin{figure}[htbp]
\centering
\begin{tabular}{cc}
\includegraphics[width=5.7 cm]{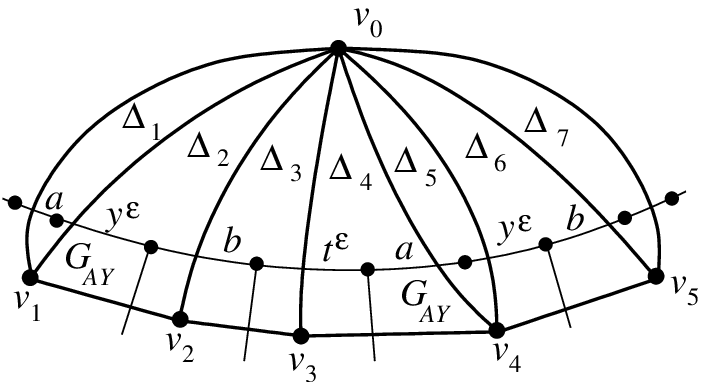}
 &
\includegraphics[width=6 cm]{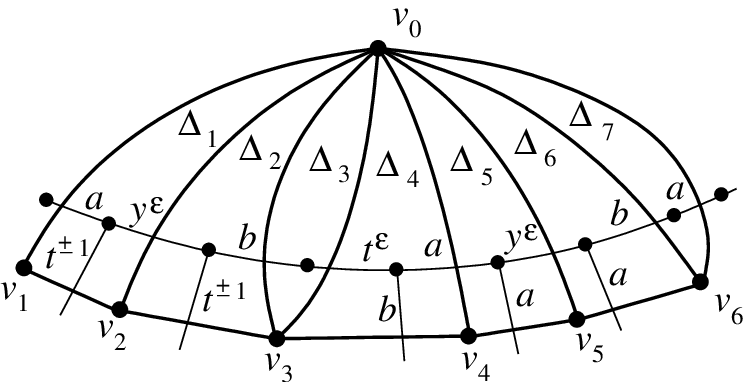}
\\
(a) & (b)
\end{tabular}
\caption{}\label{example1}
\end{figure}

\begin{proof}
Note that if $d(\Delta_i)>\pi/3$ then $i=7$. Hence,
if $|\{i\,|\,q_i=2\}|\leq 2$ then $d(S)\leq \pi/2+\pi/3-5\pi/6=0$.

Suppose that $|\{i\,|\,q_i=2\}|=3$.
If $q_i\geq 4$ for some $i$ then
$d(S)\leq 2\pi/3+\pi/2+3(-\pi/6)+(-2\pi/3)=0$.
Hence, we may assume that $q_i\leq 3$ for $1\leq i\leq 7$.
However, if $q_i=q_{i+2}=2$ for any $i\leq 4$ then $q_{i+1}\geq 4$;
moreover, if $q_7\not=2$ then this condition holds for some~$i$.
Thus $\{i\,|\,q_i=2\}=\{1,4,7\}$, $\{1,5,7\}$, $\{2,5,7\}$.
Label consecutive vertices of $S$ by $v_2,\dots,v_5$ so that
$v_2$ is adjacent to $v_1$.

\smallskip\noindent
{\bf \{1,\,4,\,7\}.}
The vertices $v_1$ and $v_3$ arise from maximal
$G_{AY}$- and $G_{AT}$-regions, respectively.
Moreover, at most two of
the five angles in
$\Delta_5$, $\Delta_6$, and $\Delta_7$ can be greater than
$\pi/3$. Then
$$
d(S)\leq \phi+(\phi+\theta-\pi)+(\theta+\psi-\pi)
+\psi+3\frac\pi3+2\frac\pi2-2\pi=2\phi+2\psi+2\theta-2\pi=0.
$$

\noindent
{\bf \{1,\,5,\,7\}.}
Since $v_1$ and $v_4$ arise from maximal
$G_{AY}$-regions and $q_2=q_3=q_4=3$,
either $v_2$ arises from a maximal $G_{AB}$-region
and $v_3$ arises from a maximal $G_{AT}$-region,
or $v_2$ arises from a maximal $G_{YB}$-region
and $v_3$ arises from a maximal $G_{YT}$-region,
see Figure~\ref{example1}(a).
In both cases
$$
d(S)\leq 5\phi+2\psi+4\theta-4\pi=
3\phi+2\theta-2\pi\leq 0.
$$

\noindent
{\bf \{2,\,5,\,7\}.}
Since $v_2$ and $v_4$ arise from maximal
$G_{BY}$- and $G_{AY}$-regions, we immediately get
$$
d(S)\leq 3\phi+3\psi+5\theta-4\pi=2\theta-\pi\leq0.
$$

\medskip

Finally, suppose that $|\{i\,|\,q_i=2\}|=4$. Then
$\{i\,|\,q_i=2\}=\{1,3,5,7\}$ and hence
$q_2\geq 4$ and $q_4\geq 4$.
Since $d(\Delta_i)\leq\phi$ for $i=1,3,5$
and $d(\Delta_7)\leq\pi/2$, we get
$d(S)\leq  3\phi+\pi/2-4\pi/3-\pi/6\leq 0$.
\end{proof}

\begin{claim}\label{claim2}
If $v_1$ arises from a maximal $G_{AT}$--region
then $d(S)\leq -\pi/3$.
\end{claim}

\begin{proof}
Since $v_1$ arises from a maximal $G_{AT}$-region,
$\Delta_1$ is not a 2-gon and, therefore, $|\{i\,|\,q_i=2\}|\leq 3$.
If $|\{i\,|\,q_i=2\}|\leq 1$ then
$d(S)\leq \pi/2+6(-\pi/6)=-\pi/2$.

Suppose that $|\{i\,|\,q_i=2\}|=2$. If
$q_i\geq 4$ for some $i$ then $d(S)\leq -\pi/2$.
Hence, we may assume that $q_i\leq 3$ for $1\leq i\leq 7$.
However, if $q_i=q_{i+2}=2$ for any $i\leq 4$ then $q_{i+1}\geq 4$;
moreover, if $q_2=2$ then $q_1\geq 4$.
This reduces us to the four cases:
$\{i\,|\,q_i=2\}=\{3,6\}$, $\{3,7\}$, $\{4,7\}$, $\{5,7\}$.
Label consecutive vertices of $S$ by $v_2,\dots,v_6$ so that
$v_2$ is adjacent to $v_1$.

\smallskip\noindent
{\bf \{3,\,6\}.}
Since $v_3$ arises from a maximal $G_{BT}$-region
and $v_5$ arises from a maximal $G_{YB}$-region,
we have
$d(S)\leq 4\psi+3\phi+5\theta-5\pi\leq -2\pi/3$.

\smallskip\noindent
{\bf \{3,\,7\}.}
If $u^\varepsilon$ is followed by $t^{p_i}$, then both
$v_3$ and $v_6$ arise from maximal $G_{BT}$-regions and thus
\begin{eqnarray*}
d(S)&=&d(\Delta_1\Delta_2\Delta_3)+d(\Delta_4\Delta_5\Delta_6)+
d(\Delta_7)\\
&\leq& (\psi+2\phi+2\theta-2\pi)
-3\frac{\pi}6+\phi\leq -\frac\pi3.
\end{eqnarray*}
So suppose that $u^\varepsilon$ is followed by~$a$
(see Figure~\ref{example1}(b)). Then
$v_6$ arises from a maximal $G_{AB}$-region.
If $v_4$ arises from a maximal $G_{AT}$-region or
$v_5$ arises from a maximal $G_{BY}$-region then
$d(S)\leq 3\psi+3\phi+6\theta-5\pi\leq -\pi/2$.
Hence, we may assume that $v_4$ arises from
a maximal $G_{AB}$-region and $v_5$ arises from
a maximal $G_{AY}$-region. It follows that ${\rm deg}(v_4)\geq 5$, that is,
the angle at $v_4$ is at most $2\pi/5$, so we have
{\setlength\arraycolsep{2pt}
\begin{eqnarray*}
d(S)&=&d(\Delta_1\Delta_2\Delta_3)+d(\Delta_4\Delta_5)+
d(\Delta_6\Delta_7)\\
&\leq& \left(\psi+2\phi+2\theta-2\pi\right)+
\left(2\phi+2\frac{2\pi}5-2\pi\right)
+\left(-\frac\pi6+\frac\pi2\right)
\leq -\frac{11\pi}{30}.
\end{eqnarray*}}

\noindent
{\bf \{4,\,7\}.}
Since both $v_1$ and $v_4$ arise from maximal $G_{AT}$-regions,
$d(\Delta_1\Delta_2\Delta_3)\leq 2\psi+2\phi+2\theta-3\pi=-\pi$.
Then
$
d(S)=d(\Delta_1\Delta_2\Delta_3)+d(\Delta_4\Delta_5\Delta_6\Delta_7)
\leq-\pi+\psi-2\pi/6+\pi/2\leq-\pi/2$.

\smallskip\noindent
{\bf \{5,\,7\}.} If $u^\varepsilon$ is followed by $t^{p_i}$, then
$q_5\geq4$, a contradiction.
Suppose $u^\varepsilon$ is followed by~$a$.
The vertex $v_5$ arises from a maximal $G_{AY}$-region
and $v_6$ arises from a maximal $G_{AB}$-region.
If $v_3$ arises from a maximal $G_{BY}$-region
or $v_4$ arises from a maximal $G_{AT}$-region then
$d(S)\leq 3\psi+3\phi+6\theta-5\pi\leq -\pi/2$. Therefore, we may assume
that $v_4$ arises from a maximal $G_{YT}$-region. It follows that
$v_3$ arises from a maximal $G_{BT}$-region and hence
${\rm deg}(v_4)\geq 5$. Thus,
$$
d(S)\leq \psi+5\phi+4\theta+2\frac{2\pi}5-5\pi\leq-\frac{11\pi}{30}.
$$

\medskip
Now suppose that $|\{i\,|\,q_i=2\}|=3$.
Then $\{i\,|\,q_i=2\}$ is one of the following:
$\{2,4,6\}$, $\{2,4,7\}$, $\{2,5,7\}$, $\{3,5,7\}$.
However, if $q_2=q_4=2$ then
$q_1\geq 4$ and $q_3\geq 4$, and so $d(S)\leq -\pi/2$.

\smallskip\noindent
{\bf \{2,\,5,\,7\}.}
Since $q_1\geq 4$, we may assume that $q_i\leq 3$ for all
$1<i\leq 7$. In particular, $q_6=3$, so
the 2-gon $\Delta_7$ comes from a maximal $G_{AB}$-region.
Since $\Delta_2$ and $\Delta_5$ come from maximal
$G_{BY}$- and $G_{AY}$-regions, respectively,
we have
$$
d(S)\leq 3\phi+4\psi+5\theta-5\pi=\psi+2\theta-2\pi\leq-\frac{2\pi}3.
$$

\noindent
{\bf \{3,\,5,\,7\}.} Then $q_4\geq 4$ and again we may assume that
$q_i\leq 3$ for all $i\not=4$ and therefore
the 2-gon $\Delta_7$ comes from a maximal $G_{AB}$-region.
Since $\Delta_3$ and $\Delta_5$ come from maximal
$G_{BT}$- and $G_{AY}$-regions, respectively,
we have
$$
d(S)\leq 6\phi+\psi+5\theta-5\pi=5\phi+4\theta-4\pi\leq-\frac{\pi}3.
$$
\end{proof}

If $w(t,u)=u^{q_1}$ then Claims~\ref{claim1} and~\ref{claim2} imply
the required contradiction that
$\sum_{\mathrm{ext. }\Delta}d(\Delta)\leq 0$.
Hence, $u$ has infinite order in~$G$.

Now suppose $w(t,u)=t^{p_1}u^{q_1}\dots t^{p_m}u^{q_m}$.
For each syllable
$t^{p_i}u^{q_i}=t^{p_i}u^\varepsilon u^{q_i-\varepsilon}$,
consider the part of the boundary corresponding to
$t^{p_i}u^\varepsilon$.
Label the first exterior region $\Delta_0$ and, as before, label the
remaining regions $\Delta_1,\dots,\Delta_7$.

If $q_0\not=2$ then $d(\Delta_0)\leq -\pi/6$
and, therefore, $d(\Delta_0S)=d(\Delta_0)+d(S)<0$.
If $\Delta_0$ is a 2-gon then $v_1$ arises from a $G_{AT}$-region.
But then, by Claim~\ref{claim2}, $d(S)\leq-\pi/3$ and, hence,
$d(\Delta_0S)\leq \psi-\pi/3\leq0$.
It follows that the sum of the $d$-values of all exterior regions
arising from any syllable $t^{p_i}u^{q_i}$ is non-positive,
and the required contradiction follows.
Hence, $t$ and $u$ generate a free product in $G$ and, since $u$
has infinite order,
$G$ contains a non-abelian free subgroup.

\bigskip

Now suppose $(\theta,\alpha,\beta)=(\pi/2,\pi/2,0)$.
Then $G$ is isomorphic to
an amalgamated free product $L*_K M$, where
$L=G_{14}*_{G_4} G_{34}$, $M=G_{12}*_{G_2} G_{23}$
and $K=G_1*G_3$.
We may assume that $|G_1|=|G_3|=2$ for otherwise $K$
(and hence $G$) contains a non-abelian free subgroup.
Similarly we may assume that $|G_2|=|G_4|=2$, so
$G_i=\langle x_i\,|\,x_i^2\rangle$ for all $i$.
Further, each $m_{i,i+1}=4$
so $G_{i,i+1}=\langle x_i,x_{i+1}\,|\,
x_i^2,x_{i+1}^2,(x_ix_{i+1})^2\rangle\cong D_4$.
Therefore, $G$ has presentation
$$
\langle x_1,x_2,x_3,x_4\, |\, x_1^2,x_2^2,x_3^2,x_4^2,(x_1 x_2)^2,
(x_2 x_3)^2,(x_3 x_4)^2,(x_4 x_1)^2\rangle.
$$
Since $G$ is a group of isometries of the Euclidean
plane, it is virtually abelian.

Thus the theorem is proved when $G$ is based on a graph with
four vertices. To complete the proof in the general case,
it remains to note that when $\G$ has five or more vertices 
it is impossible to label the edges of $\G$ so that
all four vertex subgraphs give rise to the virtually abelian
group. Therefore, one of the four vertex 
subgraph groups contains a 
non-abelian free subgroup.
Since by \cite{Cor96} subgraph groups embed,
Theorem~\ref{thm} is proved.

\section{Application}

We consider the following class of groups
which generalizes the groups defined by
periodic paired relations~\cite{Vin97}. Let $n\geq 3$,
$1\leq i,j\leq n$, $n_{ij}\geq 1$ and $1\leq t\leq n_{ij}$.
For each such $i,j,t$ let $2\leq q_i,q_{i,j;t}\leq\infty$
and suppose
$w_{i,j;t}(x_i,x_j)$ is a cyclically reduced word in $x_i$
and $x_j$. Define
$$
\Gamma=\langle x_1,\dots,x_n\,|\,
x_i^{q_i}, w_{i,j;t}(x_i,x_j)^{q_{i,j;t}}
(1\leq i,j\leq n, 1\leq t\leq n_{ij})\rangle.
$$

Each group $\Gamma$ can be realized as a Pride group by setting
$G_i=\langle x_i\,|\,x_i^{q_i}\rangle$ and $R_{ij}=\{w_{i,j;t}(x_i,x_j)^{q_{i,j;t}}\,|\,
1\leq t\leq n_{ij}\}$. For each $i,j$ define $r_{ij}={\rm min}
\{\ell_{i,j;t}q_{i,j;t}\,|\, 1\leq t\leq n_{ij}\}$, where
$\ell_{i,j;t}$ denotes the free product length 
of~$w_{i,j;t}(x_i,x_j)$. If
$1/r_{ij}+1/r_{jk}+1/r_{ik}\leq 1/2$ for all distinct $1\leq
i,j,k\leq n$ then by the Spelling Theorem for generalized triangle
groups \cite{HK06}, the Pride group $\Gamma$ is non-spherical.

\bigskip\noindent
{\bf Corollary.}
{\it
Let $\Gamma$ be as defined above with $n\geq 4$. If
$1/r_{ij}+1/r_{jk}+1/r_{ik}\leq 1/2$ for all distinct $i,j,k$, then
$\Gamma$ contains a non-abelian free subgroup unless it has
presentation
$$
\langle x_1,x_2,x_3,x_4\, |\, x_1^2,x_2^2,x_3^2,x_4^2,(x_1 x_2)^2,
(x_2 x_3)^2,(x_3 x_4)^2,(x_4 x_1)^2\rangle,
$$
in which case $\Gamma$ is virtually abelian.
}

\bigskip\noindent
{\it Acknowledgements.}
We would like to thank the referee for the careful 
reading of this paper.
The second author would like to thank the 
department of mathematics
at the Universit\'e de Provence for its hospitality during a
research visit in July 2006, when part of this work was carried out.

\end{document}